\documentstyle[12pt,epsf]{article}
\newcommand{\rr}{{\rm{I \! R}}}
\newcommand{\rn}{{\rm{I \! N}}}
\newcommand{\al}{\alpha}
\newcommand{\ds}{\displaystyle}
\newcommand{\f}{\frac}
\newcommand{\lam}{\lambda}

\newcommand{\om}{\omega}

\newcommand{\De}{\Delta}
\newcommand{\p}{\partial}
\newcommand{\ca}{{\cal A}}
\newcommand{\th}{\theta}
\newcommand{\cc}{{\rm{{\footnotesize{l}}\!\!\! C}}}

\title{\bf A dynamic p53-mdm2 model with distributed time delay}
\author{\small{M. NEAM\c TU$^{a}$\thanks{Corresponding author} , D. OPRI\c S$^{b}$, R.F. HORHAT$^c$}}
\date{}

\begin{document}
\maketitle

\begin{tabular}{cccccccc}
\scriptsize{$^{a}$Department of Economic Informatics and Statistics, Faculty of Economics,}\\
\scriptsize{West University of Timi\c soara, str. Pestalozzi, nr. 16A, 300115, Timi\c soara, Romania,}\\
\scriptsize{E-mail:mihaela.neamtu@fse.uvt.ro,}\\
\scriptsize{$^{b}$ Department of Applied
Mathematics, Faculty of Mathematics,}\\
\scriptsize{West University of Timi\c soara, Bd. V. Parvan, nr. 4, 300223, Timi\c soara, Romania,}\\
\scriptsize{E-mail: opris@math.uvt.ro}\\
\scriptsize{$^{c}$ Department of Biophysics and Medical Informatics,}\\
\scriptsize{University of Medicine and Pharmacy, Piata Eftimie Murgu, nr. 3, 300041, Timi\c soara, Romania,}\\
\scriptsize{E-mail: rhorhat@yahoo.com}\\
\end{tabular}

\bigskip

\begin{center}
\small{\bf Abstract}
\end{center}

\begin{quote} \small{ The objective of this paper is to investigate
the stability of limit cycles of a mathematical model with a
distributed delay which describes the interaction between p53 and
mdm2. Choosing the delay as a bifurcation parameter we study the
direction and stability of the bifurcating periodic solutions
using the normal form and the center manifold theorem. Some
numerical examples are finally made in order to confirm the
theoretical results.}
\end{quote}

\noindent{\small{{\it Keywords:} delay differential equation,
stability, Hopf bifurcation, P53, MDM2.

}}

\section*{\normalsize\bf 1. Introduction}

\hspace{0.6cm} The tumor suppresser gene p53 plays a key role in
oncogenesis and its anomalies are almost universal in tumoral
cells [4]. To make the reading easier we briefly present the role
of p53 gene in living cells. Normally the activity of p53 is
inhibited and the concentration of  p53 protein are kept at very
low levels. The activation occurs when there are DNA damages
[2,3]. Depending on these damages there are two outcomes: one is
the cell cycle arrest induced by a low level or a brief elevation
of p53 protein, and the other is the apoptosis induced by a high
level or a prolonged elevation of p53 protein [5]. Knowing these 2
different outcomes, it is now clear that the level of p53 protein
should be kept tight control. This control is achieved with the
help of mdm2 gene with which p53 makes a feedback loop [5, 10].
Recently it had been discovered that this loop is not quite
straightforward because there are two isoforms of the mdm2
protein, i.e. p76mdm2 and p90mdm2 which have different features
and roles [7].

This new model is based on [8, 9]. Here we achieve a smoother
modeling of the phenomenon, i.e. the interaction p53-mdm2. The
production of p53 protein is continuous, so is the binding between
p53 and the promoter of the mdm2. The difference from the previous
model [9] lies in the introduction of the integral form in the
third equation, which is the natural way of modeling a continuous
process. To have a clear picture of this aspect, we explain the
need for the introduction of the integral form, by the fact that
the synthesis of p53 protein and the binding p53-mdm2 are
continuous and there is no buffer to store p53 and then, at a
given moment, to release it, and after that all the quantity of
p53 to bind mdm2 promoter. As a matter of fact, the molecules of
p53 protein that are bound with mdm2 promoter were not synthesized
at the moment t, but they were synthesized at a previous different
moments,  they were bound at different moments.

The variables of the model are: $x_1$, $x_2$ mRNA concentrations
and $y_1$, $y_2$ the protein concentrations.

The mathematical system which describes our model is:
$$\begin{array}{l}
\vspace{0.1cm}
\dot x{}_1(t)=1-b_1x_1(t),\\
\vspace{0.1cm}
\dot y{}_1(t)=x_1(t)-(a_1+a_{12}y_2(t))y_1(t),\\
\dot x{}_2(t)=\ds\f{1}{\tau}\int_0^\tau (\alpha f(y_1(t))+(1-\alpha)f(y_1(t-s)))ds-b_2x_2(t),\\
\dot y{}_2(t)=x_2(t)-(a_2+a_{12}y_1(t))y_2(t)\end{array}\eqno(1)$$
where:  $b_1, b_2$ are the rates for mRNA degradation, $a_1, a_2,
a_{12}$ are the rates for protein degradation. The function
$f:\rr_+\rightarrow\rr$, is the Hill function, given by:
$$f(x)=\f{x^n}{a+x^n}\eqno(2)$$
with $n\in\rn^*, a>0.$ The parameters $a_1$, $a_2$, $b_1$, $b_2$,
$a_{12}$ of the model are assumed to be positive numbers less or
equal to 1, $\alpha\in[0,1]$ and $\tau\geq0.$

For $\alpha=0$ in the present model, we obtain the model from [9],
which suggests that there is an oscillatory behavior based on
using only numerical simulations.

For the study of the model (1) we consider the following initial
values:

$$x_1(0)=\bar x_1, y_1(\th)=\varphi_1(\th), \th\in[-\tau,0], x_2(0)=\bar x_{2},
y_2(0)=\bar y_2,$$ with $\bar x_1\geq0$,  $\bar x_2\geq0$, $\bar
y_2\geq0$, $\varphi_1(\theta)\geq0$, for all $\theta\in[-\tau, 0]$
and $\varphi_1$ is a differentiable function.

The paper is organized as follows. In section 2, we discuss the
local stability for the equilibrium state of system (1). We
investigate the existence of the Hopf bifurcation for system (1)
using time delay as the bifurcation parameter. In section 3, the
direction of Hopf bifurcation is analyzed by the normal form
theory and the center manifold theorem introduced by Hassard [6].
Numerical simulations for justifying the theoretical results are
illustrated in section 4. Finally, some conclusions are made and
further research directions are presented.

\section*{\normalsize\bf{2. Local stability and the existence of the Hopf bifurcation.}}

\hspace{0.6cm} The equilibrium points of system (1) are given by
the solutions of the following system of equations:

$$\begin{array}{l}
\vspace{0.1cm}
1-b_1x_1(t)=0,\\
\vspace{0.1cm}
x_1(t)-(a_1+a_{12}y_2(t))y_1(t)=0,\\
f(y_1(t))-b_2x_2(t)=0,\\
x_2(t)-(a_2+a_{12}y_1(t))y_2(t)=0.\end{array}\eqno(3)$$

Let $g:(0,\infty)\rightarrow\rr$ be the function given by:
$$g(x)=\ds\f{(a_2+a_{12}x)(1-a_1b_1x)}{b_1a_{12}x}.\eqno(4)$$ From
(3) and (4) it results that one solution of system (3) is:
$$x_{1}=\f{1}{b_1},\quad y_{1}, \quad x_{2}=\f{(a_2+a_{12}y_1)(-a_1b_1y_1+1)}{b_1a_{12}y_1}, \quad
y_{2}=\ds\f{-a_1b_1y_1+1}{b_1a_{12}y_1},\eqno(5)$$ where $y_1$ is
the solution of the equation:
$$f(x)-g(x)=0.\eqno(6)$$ Because $f(x)$ is an increasing function with $\lim\limits_{x\to \infty}
f(x)=1$ then $g(x)$ is an increasing function too, for
$x<\ds\f{1}{a_1b_1}$. It results that equation (6) has a unique
solution.

{\bf Proposition 1}. {\it If $y_{10}$ is a real solution of (6)
then the equilibrium point of system (1) is:
$$x_{10}=\ds\f{1}{b_1}, \quad y_{20}=\ds\f{1-a_1b_1y_{10}}{b_1a_{12}y_{10}},
\quad
x_{20}=\ds\f{(a_2+a_{12}y_{10})(1-a_1b_1y_{10})}{b_1a_{12}y_{10}}.\eqno(7)$$}

We consider the following translation:
$$x_1=u_1+x_{10}, y_1=u_2+y_{10}, x_2=u_3+x_{20},
y_2=u_4+y_{20}.\eqno(8)$$ With respect to (8), system (1) can be
expressed as:
$$\begin{array}{l}
\vspace{0.1cm}
\dot u{}_1(t)=-b_1u_1(t),\\
\vspace{0.1cm}
\dot u{}_2(t)=u_1(t)-(a_1+a_{12}y_{20})u_2(t)-a_{12}y_{10}u_4(t)-a_{12}u_2(t)u_4(t),\\
\dot u{}_3(t)\!=\!\ds\f{1}{\tau}\int_0^\tau (\alpha f(\!u_2(t)\!+\!y_{10}\!)\!+\!(\!1\!-\!\alpha\!)f(\!u_2(t\!-\!s\!)\!+\!y_{10}\!)\!)ds\!-\!b_2u_3(t)\!-\!b_2x_{20},\\
\dot
u{}_4(t)=-\!a_{12}y_{20}u_2(t)+u_3(t)\!-\!(a_2+a_{12}y_{10})u_4(t)\!\!-\!a_{12}u_2(t)u_4(t).\end{array}\eqno(9)$$

System (9) has $(0,0,0,0)$ as equilibrium point.

To investigate the local stability of the equilibrium state we
linearize system (9). We expand it in a Taylor series around the
origin and neglect the terms of higher order than the first order
for the functions from the right side of (9). We obtain:
$$\dot U(t)=AU(t)+\ds\f{1}{\tau}B\int_0^{\tau}U(t-s)ds,\eqno(10)$$
where
$$A\!\!=\!\!\left(\!\!\!\!\begin{array}{cccc}
\vspace{0.2cm}
-b_1 & 0 & 0 & 0\\
\vspace{0.2cm}
1 & -(a_1\!\!+\!\!a_{12}y_{20}) & 0 & -a_{12}y_{10}\\
0 & \alpha\rho_1 & -b_2 & 0\\
\vspace{0.2cm} 0 & -a_{12}y_{20} & 1 & -(a_2\!\!+\!\!a_{12}y_{10})
\end{array}\!\!\!\!\right)\!\!, B\!\!=\!\!\left(\!\!\!\!\begin{array}{cccc}
\vspace{0.2cm}
0 & 0 & 0 & 0\\
\vspace{0.2cm}
0 & 0 & 0 & 0\\
\vspace{0.2cm}
0 & (\!1\!-\!\alpha\!)\rho_1 & 0 & 0\\
\vspace{0.2cm} 0 & 0 & 0 & 0\end{array}\!\!\!\!\right)\eqno(11)$$
with $\rho_1=f'(y_{10})$,
$U(t)=(u_1(t)$,$u_2(t)$,$u_3(t)$,$u_4(t))^T$, \vspace{0.2cm}
$\int_0^{\tau}U(t-s)ds=(\int_0^{\tau}u_1(t-s)ds$,$\int_0^{\tau}u_2(t-s)ds$,$\int_0^{\tau}u_3(t-s)ds$,
$\int_0^{\tau}u_4(t-s)ds)^T.$

\smallskip

The characteristic equation corresponding to system (10) is
$det(\lambda I-A-(\ds\f{1}{\tau}\int_0^{\tau}e^{-\lambda
s}ds)B)=0$ which leads to:
$$(\lam+b_1)(\lam^3+b\lam^2+c\lam+d+\ds\f{h}{\tau}\int _0^\tau e^{-\lam s}ds)=0,\eqno(12)$$
where
$$\begin{array}{l}
b=a_1+a_2+b_2+a_{12}(y_{20}+y_{10}),\\
c=b_2(a_1+a_2)+b_2a_{12}(y_{20}+y_{10})+a_1a_2+a_{12}(a_1y_{10}+a_2y_{20}),\\
d=b_2a_1a_2+a_{12}b_2(y_{20}a_2+a_{1}y_{10})+\alpha a_{12}y_{10}\rho_1,\\
h=(1-\alpha)a_{12}y_{10}\rho_1 .\end{array}\eqno(13)$$

The equilibrium point $X^*=(x_{10}, y_{10}, x_{20}, y_{20})^T$ is
locally asymptotically stable if and only if all eigenvalues of
(12) have negative real parts.

Because $b_1>0$ we will analyze the function:
$$\De (\lam,
\tau)=\lam^3+b\lam^2+c\lam+d+\ds\f{h}{\tau}\int_0^{\tau}e^{-\lam
s}ds,\eqno(14)$$ with $\lam\in\rr$.

We are going to show that the equilibrium point $X^*$ undergoes a
Hopf bifurcation. In this sense, we look for the existence of the
purely imaginary roots of $\De (\lam, \tau)=0$. First, we verify
if $X^*$ is locally asymptotically stable when $\tau=0$. In this
case, the equation $\De (\lam, \tau)=0$ becomes:
$$\lam^3+b\lam^2+c\lam+d+h=0.\eqno(15)$$
Because the coefficients of equation (15) are positive then
according to the Routh-Hurwitz criterion we have:

{\bf Proposition 2.} {\it When there is no delay, the equilibrium
point $X^*$ of system (1) is locally asymptotically stable if and
only if $$cb>d+h,$$ where $c,b,d,h$ are given by (13).}

We are looking for the values $\tau_0$ so that the equilibrium
point $X^*$ changes from local asymptotic stability to instability
or vice versa. This is specific for the characteristic equation
with pure imaginary solutions. Let $\lam=\pm i\om$ be these
solutions. We assume $\om>0$. It is sufficient to look for
$\lam=i\om$ root of $\De (\lam, \tau)=0$. Separating real and
imaginary parts of $\De (i\om, \tau)=0$ we obtain:
$$sin(\om\tau)=\ds\f{\tau\om(b\om^2-d)}{h}, \quad
cos(\om\tau)=1-\ds\f{\tau\om^2(c-\om^2)}{h}.\eqno(16)$$

A solution of (16) is a pair $(\om,\tau)$ so that
$\om\tau\in[0,2\pi]$, $\om\in(0,\sqrt c]$ and $\tau=g_1(\om)$,
where $g_1:[0, \sqrt c]\longrightarrow\rr_+$ is given by:
$$g_1(x)=\ds\f{2h(c-x^2)}{(bx^2-d)^2+x^2(c-x^2)^2}.\eqno(17)$$
From (17) it results that $g_1(x)\in[0,\ds\f{2hc}{d^2}]$ and for
all $\tau\in[0,\ds\f{2hc}{d^2}]$ there is $\om\in[0,\sqrt c]$ so
that $\tau=g(\om)$.

We have:

{\bf Proposition 3.} {\it For $\tau\in[0,\ds\f{2hc}{d^2}]$ there
is $\om\in[0,\sqrt c]$, so that $(\om, \tau)$ is solution of
(16).}

According to the above proposition, system (16) has not unique
solution.

In order to show that $X^*$ undergoes a Hopf bifurcation for
$\tau=\tau_0$ we have to prove that $\pm i\om_0$ are simple
eigenvalues of $\De (\lam,\tau_0)$ and satisfy the transversality
condition $\ds\f{dRe(\lam)}{d\tau}|_{\tau=\tau_0}\neq0.$

From (14), it results that:
$$\De_\lam(\lam,\tau)=\ds\f{\p
\De}{\p\lam}(\lam,\tau)=3\lam^2+2b\lam+c-\ds\f{h}{\lam^2\tau}(1-(1+\lam\tau)e^{-\lam\tau}).\eqno(18)$$
From (18), for $(\om,\tau)$ a solution of (16) and $\tau=g(\om)$
it results that:
$$\begin{array}{l}
M_1=Re(\De_{\lam}(i\om,g_1(\om))=-4\om^2+2c-\tau(b\om^2-d)\\
M_2=Im(\De_{\lam}(i\om,g_1(\om))=\ds\f{3b\om^2-d-h}{\om}+\tau\om(c-\om^2).\end{array}$$
Then, we have:
$$\begin{array}{lll}
M_1^2\!\!+\!\!M_2^2 & \!\!=\!\! &
(2c\!-\!4\om^2)^2\!+\!\tau^2(b\om^2\!-\!d)^2\!+\!\ds\f{(3b\om^2\!-\!d\!-\!h)^2}{\om^2}\!+\! \tau^2\om^2(c\!-\!\om^2)^2\!+\\
 & \!\!+\!\! &
\!2\tau((3b\om^2\!-\!d\!-\!h)(c\!-\!\om^2)\!-\!(b\om^2\!-\!d)(2c\!-\!4\om^2)).\end{array}\eqno(19)$$

From (19), it results that:

{\bf Proposition 4.} If $bc>d+h$, $d>h$, $bc>3d-h$, then
$\lambda=i \omega_0$ is a simple root for the equation
$\Delta(\lambda,\tau_0)=0.$

Now, we consider a branch of eigenvalues
$\lam(\tau)=\nu(\tau)+i\om(\tau)$ of (14) so that $\nu(\tau_0)=0$
and $\om(\tau_0)=\om_0$, where $(\om_0, \tau_0)$ is a solution of
(16). Differentiating equation $\De(\lam(\tau), \tau)=0$ with
respect ro $\tau$, we obtain:
$$\lam'(\tau)=\ds\f{d\lam}{d\tau}|_{\tau=\tau_0,
\om=\om_0}=-\ds\f{\De _\tau(i\om_0, \tau_0)}{\De_\lam (i\om_0,
\tau_0)}=M(\om_0, \tau_0)+iN(\om_0,\tau_0)$$ where
$$M(\om_0,\tau_0)=-\ds\f{M_1N_1+M_2N_2}{M_1^2+M_2^2}\eqno(20)$$ with
$$\begin{array}{l}
N_1=\ds\f{1}{\tau}[h-(b\om^2-d)-\tau\om^2(c-\om^2)]\\
N_2=\ds\f{1}{\tau}[\om(c-\om^2)-\tau\om(b\om^2-d)]\end{array}$$
and
$$N=\ds\f{M_1N_2-M_2N_1}{M_1^2+M_2^2}.\eqno(21)$$
By direct calculation $M\neq0$.

We can conclude that when $\tau_0\in[0,\ds\f{2hc}{d^2}]$ the
characteristic equation $\De (\lam, \tau_0)=0$ has a unique pair
of purely imaginary simple eigenvalues satisfied $\ds\f{dRe
\lam(\tau)}{d\tau}|_{\tau=\tau_0}\neq0$ and the real roots are
negative. Consequently, a Hopf bifurcation occurs at $X^*$ when
$\tau=\tau_0$. Moreover, applying Rouche's theorem we can verify
that every eigenvalue of $\De(\lam,\tau)=0$ with $\tau<\tau_0$ has
negative real part. It follows that $X^*$ is locally
asymptotically stable for $0\leq\tau<\tau_0$. These results are
summed up in the following theorem:

{\bf Theorem 1.} {\it Assume that $cb>d+h$. Then there exists
values $\tau_0\in[0,\ds\f{2hc}{d^2}]$ of time delay so that the
equilibrium point $X^*$ is locally asympto\-ti\-ca\-lly stable
when $\tau\in[0,\tau_0)$ and becomes unstable when $\tau=\tau_0$
throughout a Hopf bifurcation. In particular, the periodic
solutions appear for system (1) when $\tau=\tau_0$.}

\section*{{\normalsize\bf 3. Direction and stability of the Hopf bifurcation}}

\hspace{0.6cm} In the previous section, we obtain some conditions
with guarantee that system (1) undergoes Hopf bifurcation at
$\tau=\tau_0$.

In this section, we study the direction, the stability and the
period of the bifurcating periodic solutions. The used method is
based on the normal form theory and the center manifold theorem
introduced by Hassard [6].

For an interval $I\subseteq \rr$, we define the space of
continuous functions as $C(I,K)=\{f:I\rightarrow K, f$
continuous$\},$ where $K=\rr^4$ or $\cc^4$. When $I=[-\tau,0]$,
$\tau=\tau_0+\mu$, $\mu>0$ sufficiently small,  we set
$C_\mu=C([-\tau,0],K)$. Expanding the functions from the right
side of system (9) in Taylor series around $(0,0,0,0)^T$ it
results that:
$$\dot X(t)\!\!=\!\!AX(t)\!\!+\!\!\f{1}{\tau}B\int_0^\tau X(t-s)ds\!\!+
\!\!F(X(t), \int_0^\tau (\alpha X(t)\!\!+\!\!(1\!\!-\!\!\alpha)X(t\!\!-\!\!s))ds)\eqno(22)$$ where

$X(t)=(u_1(t),u_2(t),u_3(t), u_4(t))^T$, $$\int_0^\tau
X(t-s)ds=(\int_0^\tau u_1(t-s)ds,\int_0^\tau u_2(t-s)ds,
\int_0^\tau u_3(t-s)ds, \int_0^\tau u_4(t-s)ds)^T,$$
$$\begin{array}{l}
F(X(t), \int_0^\tau (\alpha X(t)+(1-\alpha)X(t-s))ds)=(0, \!F^2(u_2(t)\!, \!u_4(t))\!,\\
\vspace{0.2 cm} \!F^3(\int_0^\tau (\alpha
u_2(t)+(1-\alpha)u_2(t-s))ds)\!,\!F^4(u_2(t)\!, u_4(t)))^T,
\end{array}\eqno(23)$$
 \vspace{0.2 cm}
$$\begin{array}{l}
F^2(u_2(t),u_4(t))\!\!=\!\!-a_{12}u_2(t)u_4(t),\\
\vspace{0.2 cm}
F^3(\int_0^\tau (\alpha
u_2(t)\!+\!(1\!-\!\alpha)u_2(t\!-\!s))ds)\!\!=\!\!\ds\f{1}{2\tau}\rho_2\int_0^\tau
(\alpha
u_2(t)\!+\!(1\!-\!\alpha)u_2(t\!-\!s))^2ds\!+\!\\
\vspace{0.2 cm}
+\ds\f{1}{6\tau}\rho_3\int_0^\tau (\alpha u_2(t)+(1-\alpha)u_2(t-s))^3ds,\\
\vspace{0.2 cm}
 F^4(u_2(t),u_4(t))=-a_{12}u_2(t)u_4(t),\end{array}$$ $\rho_2=f''(y_{10})$,
 $\rho_3=f'''(y_{10})$ and A,B are given by (11).

For $\Phi\in C_\mu$ with $K=\cc^4$ we define a linear operator:
$$L_\mu(\Phi)=A\Phi(0)-\f{1}{\tau_0}B\int_{-\tau_0}^0\Phi(s)ds$$
and a nonlinear operator:
$$F_\mu(\Phi)=(0, F^2(\Phi_2(0), \Phi_4(0)), F^3(\int_{-\tau}^0\Phi_2(s)ds),
F^4(\Phi_2(0),\Phi_4(0)))^T.$$

For $\Phi\in C^1([-\tau_0, 0], \cc^{4})$ we define:
$$\ca(\mu)\Phi(\th)=\left\{\begin{array}{ll} \vspace{0.2cm}
\ds\f{d\Phi(\th)}{d\th}, & \th\in[-\tau_0,0)\\
A\Phi(0)-\ds\f{1}{\tau_0}B\int_{-\tau_0}^0\Phi(s)ds, &
\th=0,\end{array}\right.$$
$$R(\mu)\Phi(\th)=\left\{\begin{array}{ll} \vspace{0.2cm}
0, & \th\in[-\tau_0,0)\\
F_\mu(\Phi), & \th=0\end{array}\right.$$ and for $\Psi\in
C^1([0,\tau_0], \cc^{*4})$, we define the adjoint operator $\ca^*$
of $\ca$ by:
$$\ca^*\Psi(s)=\left\{\begin{array}{ll} \vspace{0.2cm}
-\ds\f{d\Psi(s)}{ds}, & s\in(0, \tau_0]\\
A\Psi(0)+\ds\f{1}{\tau_0}(\int_0^{\tau_0}\Psi(\th)d\th) B, &
s=0.\end{array}\right.$$

Then, we can rewrite (22) in the following vector form:
$$\dot X_t=A(\mu)X_t+R(\mu)X_t\eqno(24)$$
where $X_t=X(t+\th)$ for $\th\in[-\tau_0, 0]$. For $\Phi\in
C([-\tau_0, 0], \cc^{*4})$ and $\Psi\in C([0,\tau_0], \cc^{*4})$
we define the following bilinear form:
$$<\Psi(s),
\Phi(\th)>=\bar
\Psi(0)\Phi(0)-\int_{-\tau_0}^0\int_{\xi=0}^\theta\bar\Psi(\xi-\theta)B(\ds\f{1}{\tau_0}\int_0^\xi\Phi(\xi
')d\xi ')d\xi d\th,$$ $s\in[0,\tau_0]$, $\theta\in[-\tau_0,0]$.

Then, it can verified that $\ca^*$ and $\ca(0)$ are adjoint
operators with respect to this bilinear form.

In the light of the obtained results in the last section, we
assume that $\pm i\om_0$ are eigenvalues of $\ca(0)$. Thus, they
are also eigenvalues of $\ca^*$. We can easily obtain:
$$\Phi(\th)=ve^{\lam_1\th},\quad \th\in[-\tau_0, 0]\eqno(25)$$
where $v=(v_1, v_2, v_3, v_4)^T$,
$$\begin{array}{l}v_1=0,
v_2=a_{12}y_{10},
v_3=a_{12}^2y_{10}y_{20}-(\lam_1+a_1+a_{12}y_{20})(\lam_1+a_2+a_{12}y_{10}),\\
v_4=-(\lam_1+a_1+a_{12}y_{20})\end{array}$$ is the eigenvector of
$\ca(0)$ corresponding to $\lam_1=i\om_0$ and
$$\Psi(s)=we^{\lam_1s},\quad s\in[0,\tau_0]$$ where
$w=(w_1, w_2, w_3, w_4)$,
$$w_1\!=\!\ds\f{1}{\bar\eta}, w_2\!=\!\ds\f{d_2}{\bar\eta},
w_3\!=\!\ds\f{d_3}{\bar\eta}, w_4\!=\!\ds\f{d_4}{\bar\eta},$$
$$\begin{array}{l}
d_2=b_1+\lam_1,
d_3=-\ds\f{a_{12}y_{10}(b_1+\lam_1)}{(\lam_1+a_2+a_{12}y_{10})(b_2+\lam_1)},
d_4=-\ds\f{a_{12}y_{10}(b_1+\lam_1)}{\lam_1+a_2+a_{12}y_{10}}\\
\vspace{0.2cm}\eta=v_2\bar d_2+v_3\bar d_3+v_4\bar d_4-\bar
d_3v_2\ds\f{(1-\alpha)\rho_1}{\tau_0\lam_1^3}(\tau_0\lam_1-2+2e^{-\lam_1\tau_0}+\lam_1\tau_0e^{-\lam_1\tau_0})
\end{array}$$ is the eigenvector of $\ca^*$ corresponding to
$\lam_2=-i\om_0$.

We can verify that: $<\Psi(s), \Phi(s)>=1$, $<\Psi(s),
\bar\Phi(s)>=<\bar\Psi(s), \Phi(s)>=0$, $<\bar\Psi(s),
\bar\Phi(s)>=1.$

Using the approach in [1], we next compute the coordinates to
describe the center manifold $\Omega_0$ at $\mu=0$. Let
$X_t=X(t+\th), \th\in[-\tau_0,0]$, be the solution of equation
(24) when $\mu=0$ and
$$z(t)=<\Psi, X_t>,
\quad w(t,\th)=X_t(\th)-2Re\{z(t)\Phi(\th)\}.$$

On the center manifold $\Omega_0$, we have:
$$w(t,\th)=w(z(t), \bar z(t), \th)$$ where
$$w(z,\bar z, \th)=w_{20}(\th)\ds\f{z^2}{2}+w_{11}(\th)z\bar
z+w_{02}(\th)\ds\f{\bar z^2}{2}+w_{30}(\th)\ds\f{z^3}{6}+\dots$$
in which $z$ and $\bar z$ are local coordinates for the center
manifold $\Omega_0$ in the direction of $\Psi$ and $\bar\Psi$ and
$w_{02}(\th)=\bar w_{20}(\th)$.

For solution $X_t$ of equation (24), as long as $\mu=0$, we have:
$$\dot z(t)=\lam_1z(t)+g(z, \bar z)\eqno(26)$$ where
$$\begin{array}{ll}
g(z, \bar z)& =\bar\Psi(0)F(w(z,\bar z, 0)+Re(z\Phi(0)))=\\
& =g_{20}\ds\f{z^2}{2}+g_{11}z\bar z+g_{02}\ds\f{\bar
z^2}{2}+g_{21}\ds\f{z^2\bar z}{2}+\dots\end{array}\eqno(27)$$ From
(23), (26) and (27) we obtain:
$$\begin{array}{l}g_{20}=F^2_{20}\bar w_2+F^3_{20}\bar w_3+F_{20}^4\bar w_4,
g_{11}=F^2_{11}\bar w_2+F^3_{11}\bar w_3+F_{11}^4\bar w_4,\\
g_{02}=F^2_{02}\bar w_2+F^3_{02}\bar w_3+F_{02}^4\bar
w_4,\end{array}\eqno(28)$$ where
$$\begin{array}{l}F_{20}^2=F_{20}^4=-2a_{12}v_2v_4, F_{11}^2=F_{11}^4=-a_{12}(v_2\bar v_4+\bar
v_2v_4), \\
F_{02}^2=F_{02}^4=-2a_{12}\bar v_2\bar v_4 \\
F_{20}^3=\ds\f{\rho_2v_2^2}{2\tau_0\lam_1}(2\al^2\tau_0\lam_1-4\al(1-\al)(e^{-\lam_1\tau_0}-1)-(1-\al)^2(e^{-2\lam_1\tau_0}-1))\\
\vspace{0.3 cm}
F_{11}^3=\ds\f{\rho_2v_2\bar v_2}{\tau_0\lam_1\lam_2}(\lam_1\lam_2(\al^2+(1-\al)^2)\tau_0-\al(1-\al)(\lam_2e^{-\lam_1\tau_0}+\lam_1e^{-\lam_2\tau_0}))\\
\vspace{0.3 cm}
 F_{02}^3=\ds\f{\rho_2\bar
v_2^2}{2\tau_0\lam_2}(2\al^2\tau_0\lam_2-4\al(1-\al)(e^{-\lam_2\tau_0}-1)-(1-\al)^2(e^{-2\lam_2\tau_0}-1))
\end{array}$$ and
$$g_{21}=F_{21}^2\bar w_2+F_{21}^3\bar w_3+F_{21}^4\bar w_4\eqno(29)$$ where
$$\begin{array}{l}
F_{21}^2=F_{21}^4=-a_{12}\bar v_2w_{20}^4(0)-2a_{12}v_2w_{11}^4(0)-a_{12}\bar v_4w_{20}^2(0)-2a_{12}v_4w_{11}^2(0)\\
\vspace{0.1 cm}
F_{21}^3=\ds\f{\rho_2}{\tau_0}(2v_2(\al^2\tau_0w_{11}^2(0)+\al(1-\al)k_1-\ds\f{\al(1-\al)w_{11}^2(0)}{\lam_1}(e^{-\lam_1\tau_0}-1)+\\
\vspace{0.1 cm} +(1-\al)^2k_2)+2\bar
v_2(\al^2\tau_0w_{20}^2(0)-\ds\f{\al(1-\al)}{\lam_2}w_{20}^2(0)(e^{-\lam_1\tau_0}-1)+\\
\vspace{0.1 cm} +\al(1-\al)k_3+(1-\al)^2
k_4))+\\
\vspace{0.1 cm}
 +\ds\f{\rho_3}{\tau_0}v_2^2\bar
v_2^2(\al^3\tau_0\!\!-\!\!\ds\f{(1\!\!-\!\!\al)^2\al}{2\lam_1}
(e^{-2\lam_1\tau_0}\!\!-\!\!1)\!\!-\!\!\ds\f{(1\!\!-\!\!\al)\al^2}{\lam_2}(e^{-\lam_2\tau_0}\!\!-\!\!1)\!\!+\!\!(1\!\!-\!\!\al)^3\tau_0)

\end{array}$$ with

$$k_1=\int_0^{\tau_0}w_{11}^2(-s)ds,
k_2=\int_0^{\tau_0}e^{-\lam_1s}w_{11}^2(-s)ds,$$
$$k_3=\int_0^{\tau_0}w_{20}^2(-s)ds,
k_4=\int_0^{\tau_0}e^{-\lam_2s}w_{20}^2(-s)ds$$

$$\begin{array}{l}
w_{20}^2(-s)=-\ds\f{g_{20}}{\lam_1}v_2e^{-\lam_1s}-\ds\f{\bar
g_{02}}{3\lam_1}\bar v_2e^{-\lam_2s}+E_2^2e^{-2\lam_1s}\\
w_{11}^2(-s)=\ds\f{g_{11}}{\lam_1}v_2e^{-\lam_1s}-\ds\f{\bar
g_{11}}{\lam_1}\bar v_2e^{-\lam_2s}+E_1^2\\
w_{20}^2(0)=-\ds\f{g_{20}}{\lam_1}v_2-\ds\f{\bar
g_{02}}{3\lam_1}\bar v_2+E_2^2\\
w_{11}^2(0)=\ds\f{g_{11}}{\lam_1}v_2-\ds\f{\bar
g_{11}}{\lam_1}\bar v_2+E_1^2\\
w_{20}^4(0)=-\ds\f{g_{20}}{\lam_1}v_4-\ds\f{\bar
g_{20}}{3\lam_1}\bar v_4+E_2^4\\
w_{11}^4(0)=\ds\f{g_{11}}{\lam_1}v_4-\ds\f{\bar
g_{11}}{\lam_1}\bar v_4+E_1^4,\end{array}\eqno(30)$$
$s\in[0,\tau_0]$, $E_1^2, E_1^4$ respectively $E_2^2, E_2^4$ are
the components of the vectors:
$$\begin{array}{l}
E_2=-(A-\ds\f{1}{2\lam_1\tau_0}(e^{-2\lam_1\tau_0}-1)B-2\lam_1I)^{-1}F_{20}\\
E_1=-(A+B)^{-1}F_{11},\end{array}$$ where $F_{20}=(0, F_{20}^2,
F_{20}^3, F_{20}^4)^T$, $F_{11}=(0, F_{11}^2, F_{11}^3,
F_{11}^4)^T$.



Based on the above analysis and calculation, we can see that each
$g_{ij}$ in (28), (29) is determined by the parameters and delay
from system (1). Thus, we can explicitly compute the following
quantities:
$$\begin{array}{l}
C_1(0)=\ds\f{i}{2\om_0}(g_{20}g_{11}-2|g_{11}|^2-\ds\f{1}{3}|g_{02}|^2)+\ds\f{g_{21}}{2}\\
\vspace{0.2cm} \mu_2=-\ds\f{Re(C_1(0))}{M(\om_0,\tau_0)},
T_2=-\ds\f{Im(C_1(0))+\mu_2N(\om_0,\tau_0)}{\om_0},
\beta_2=2Re(C_1(0)),\end{array}\eqno(31)$$ where
$M(\om_0,\tau_0)$, $N(\om_0,\tau_0)$ are given by (20) and (21).

In summary, this leads to the following result:

\vspace{2mm} {\bf Proposition 5.} {\it In formulas (31), $\mu_2$
determines the direction of the Hopf bifurcation: if $\mu_2>0
(<0)$, then the Hopf bifurcation is supercritical (subcritical)
and the bifurcating periodic solutions exist for $\tau>\tau_0
(<\tau_0)$; $\beta_2$ determines the stability of the bifurcating
periodic solutions: the solutions are orbitally stable (unstable)
if $\beta_2<0 (>0)$; and $T_2$ determines the period of the
bifurcating periodic solutions: the period increases (decreases)
if $T_2>0 (<0)$.}

\section*{\normalsize\bf 4. Numerical examples.}

\hspace{0.6cm}For the numerical simulations we use Maple 9.5. In
this section, we consider system (1) with $a_1=a_2=0.13$,
$a_{12}=0.02$, $a_{12}=0.06$, $b_1=0.2$, $b_2=0.4$, $a=4$,
$\al=0.2$, $n=3$. Waveform plot  are obtained by the formula:
$$X(t+\th)\!=\! z(t)\Phi(\th)\!+\!\bar
z(t)\bar\Phi(\th)\!+\!\ds\f{1}{2}w_{20}(\th)z^2(t)+w_{11}(\th)z(t)\bar
z(t)\!+\!\ds\f{1}{2}w_{02}(\th)\bar z(t)^2+X_0,$$ where $z(t)$ is
the solution of (26), $\Phi(\th)$ is given by (25), $w_{20}(\th),
w_{11}(\th), w_{02}(\th)$ are given by (30) and $X_0=(x_{10},
y_{10}, x_{20}, y_{20})^T$ is the equilibrium state.

We obtain: $x_{10}= 5$, $y_{10}= 21.03417191$, $y_{20}\!=\!
1.795140515$, $x_{20}\!=\! 2.498925919$, $\mu_2\!=\!
-0.2101567953$, $\beta_2\!=\! -0.3029980114$, $T_2\!=\!
0.1148699183$, $\omega\!=\! 0.1$, $\tau\!=\! 0.1001651263$. Then
the Hopf bifurcation is subcritical, the solutions are orbitally
stable and the period of the solution is increasing. The wave
plots are displayed in fig1 and fig2:

\begin{center}
{\small \begin{tabular}{c|c} \hline Fig.1. $(t,y_1(t))$&Fig.2.
$(t,y_2(t))$\\&\\
 \cline{1-2} \epsfxsize=5cm

\epsfysize=5cm

\epsffile{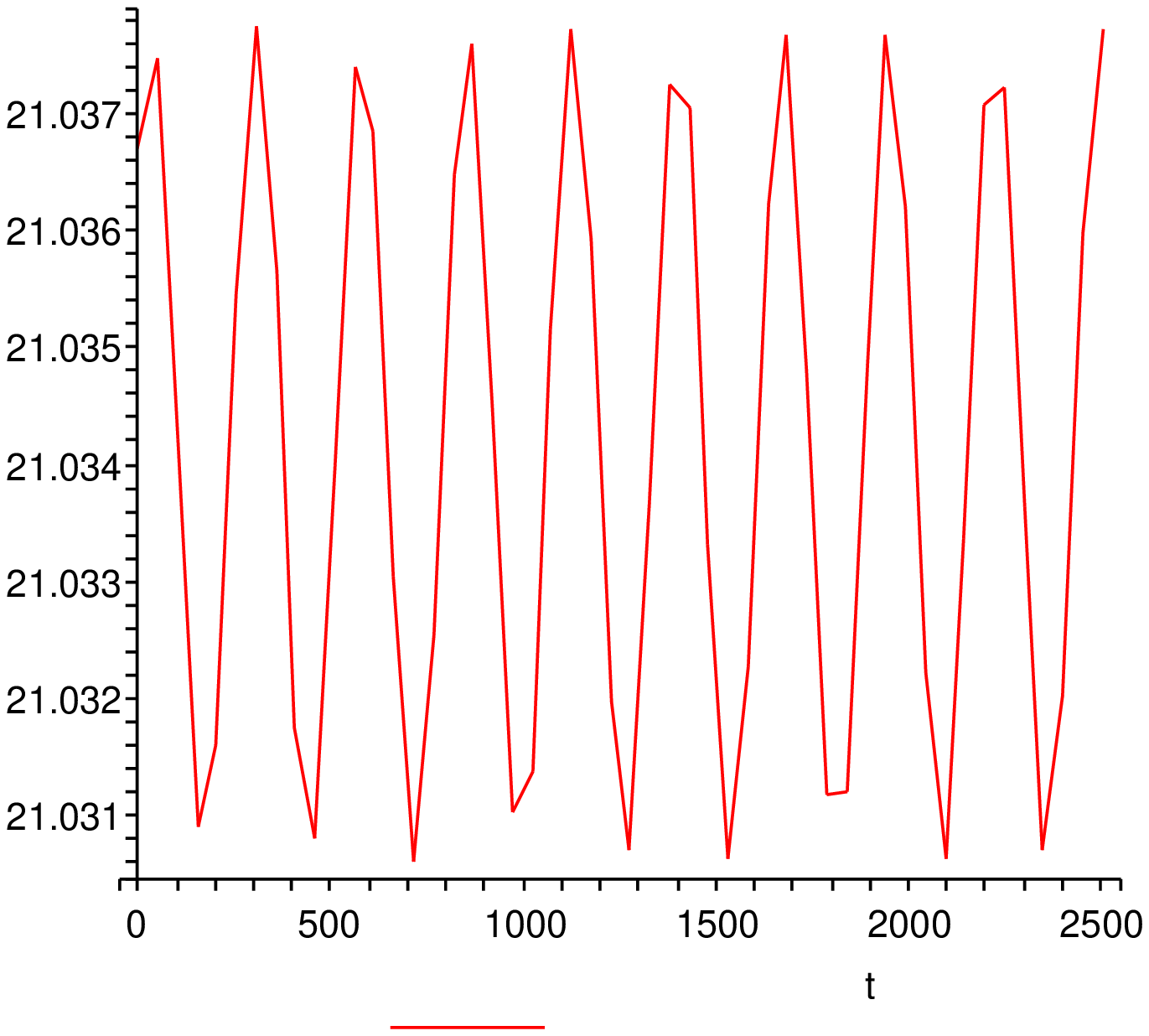} &

\epsfxsize=6cm

\epsfysize=5cm

\epsffile{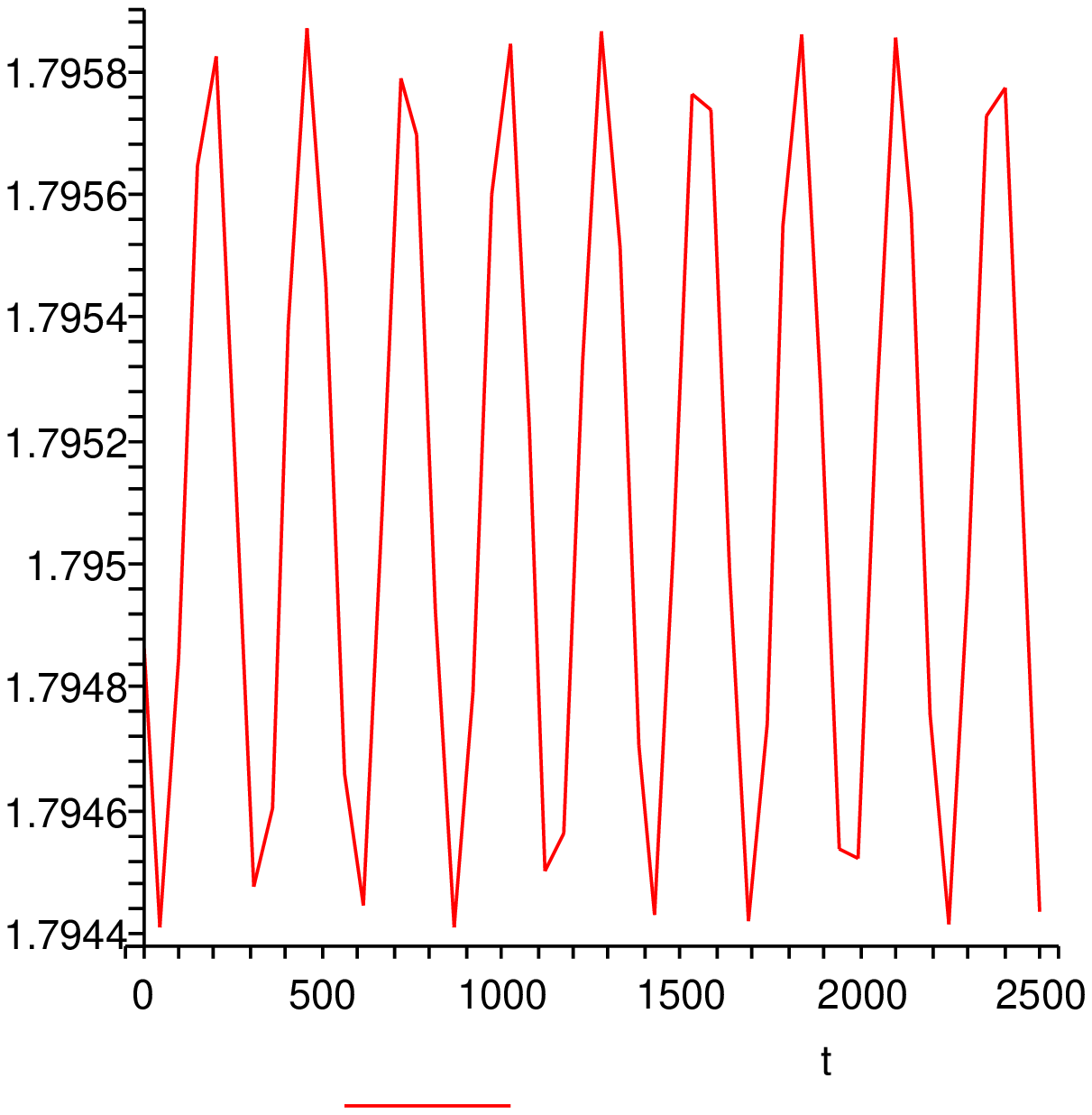}

\\
 \hline
\end{tabular}}
\end{center}

\medskip

\section*{\normalsize\bf 5. Conclusions.}

\hspace{0.6cm} As in our previous model [9], we obtain an
oscillatory behavior similar to that observed experimentally [3].
The conclusion is not surprising, but is useful as this model
provides a more accurate approach of the interaction p53-mdm2. We
can conclude that the transformation made by us to the continuous
model with distributed time of the interaction p53-mdm2, which
actually is more real, did not alter the behavior of the system.

\end{document}